\documentclass{IEEEtran4PSCC}
\ifCLASSINFOpdf
  \usepackage[pdftex]{graphicx}
\else

   \usepackage[dvips]{graphicx}

\fi

\usepackage[cmex10]{amsmath}
\usepackage{amssymb}
\usepackage{url}
\usepackage{booktabs}
\usepackage{tikz}
\usetikzlibrary{arrows.meta,positioning,fit,shapes.multipart,calc}
\usetikzlibrary{decorations.pathreplacing}

\makeatletter
\let\old@ps@headings\ps@headings
\let\old@ps@IEEEtitlepagestyle\ps@IEEEtitlepagestyle
\def\psccfooter#1{%
    \def\ps@headings{%
        \old@ps@headings%
        \def\@oddfoot{\strut\hfill#1\hfill\strut}%
        \def\@evenfoot{\strut\hfill#1\hfill\strut}%
    }%
    \def\ps@IEEEtitlepagestyle{%
        \old@ps@IEEEtitlepagestyle%
        \def\@oddfoot{\strut\hfill#1\hfill\strut}%
        \def\@evenfoot{\strut\hfill#1\hfill\strut}%
    }%
    \ps@headings%
}
\makeatother

\usepackage[none]{hyphenat} 
\usepackage{microtype}      
\begin{document}

\title{Explicit Ensemble Learning Surrogate for Joint Chance-Constrained Optimal Power Flow}

\author{
\IEEEauthorblockN{Amir Bahador Javadi}
\IEEEauthorblockA{Department of Electrical and Computer Engineering \\
Louisiana State University (LSU)\\
Baton Rouge, LA, USA\\
\url{ajavad2@lsu.edu}}
\and
\IEEEauthorblockN{Amin Kargarian}
\IEEEauthorblockA{Department of Electrical and Computer Engineering \\
Louisiana State University (LSU)\\
Baton Rouge, LA, USA\\
\url{kargarian@lsu.edu}}}

\maketitle

\begin{abstract}
The increasing penetration of renewable generation introduces uncertainty into power systems, challenging traditional deterministic optimization methods. Chance-constrained optimization offers an approach to balancing cost and risk; however, incorporating joint chance constraints introduces computational challenges. This paper presents an ensemble support vector machine surrogate for joint chance constraint optimal power flow, where multiple linear classifiers are trained on simulated optimal power flow data and embedded as tractable hyperplane constraints via Big--M reformulations. The surrogate yields a polyhedral approximation of probabilistic line flow limits that preserves interpretability and scalability. Numerical experiments on the IEEE 118-bus system show that the proposed method achieves near-optimal costs with a negligible average error of $0.03\%$. These results demonstrate the promise of ensemble surrogates as efficient and transparent tools for risk-aware optimization of power systems.
\end{abstract}


\thanksto{\noindent This work is supported by the National Science Foundation under Grant Number 1944752.}

\section{Introduction}

The growing integration of renewable energy sources (RES), particularly wind and solar generation, has profoundly reshaped modern power system operations. Renewable generation offers substantial environmental and economic benefits. However, its inherent variability and uncertainty introduce significant challenges to system reliability, economic efficiency, and secure operation \cite{morales2013integrating}. Traditional deterministic optimization methods are often inadequate in this context, as they fail to capture the stochastic nature of RES and demand fluctuations. 

Chance-constrained programming \cite{charnes1959chance} has emerged as an approach for explicitly incorporating uncertainty into optimization, allowing for controlled levels of risk. In the context of power systems, chance-constrained optimal power flow (CC-OPF) ensures that operational constraints, such as generator and line flow limits, are satisfied with a prescribed probability \cite{bienstock2014chance}. This enables operators to balance risk and cost rather than relying on worst-case assumptions that may lead to overly conservative and uneconomic dispatch solutions. However, the practical application of CC-OPF remains challenging. Exact reformulations are typically tractable only under assumptions, such as Gaussian uncertainty and linearized system models \cite{venzke2020convex,dall2017chance}. Joint chance constraints (JCCs), which enforce simultaneous satisfaction of multiple limits, are inherently non-linear and often NP-hard \cite{roald2015security}. Approximate techniques, such as Boole’s inequality or scenario approximation, mitigate complexity but often introduce excessive conservatism or computational burdens \cite{ding2022distributionally}.

To overcome these challenges, recent research has focused on surrogate models and machine learning (ML) methods that approximate the feasible region induced by network constraints under uncertainty. By learning from data, surrogates can replace intractable chance constraints with tractable approximations, reducing computational complexity while maintaining acceptable accuracy. Examples include deep quantile regression for estimating line flow distributions and constructing probabilistic constraints \cite{9956906}, convex relaxations that provide tractable reformulations of probabilistic constraints \cite{venzke2020convex}, and classification-based surrogates such as support vector machines (SVMs) \cite{vapnik1995support}, which separate feasible from infeasible operating points with scalable performance. Embedding SVM-derived constraints into the optimal power flow (OPF) formulation allows operators to enforce network security constraints with high fidelity while retaining computational efficiency, and ensemble strategies further reduce the risk of misclassification \cite{javadi2025learning}.

Although these black-box ML approaches are powerful, they often lack interpretability and fail to reveal explicit structural relationships among system variables. Symbolic regression (SR) has recently gained attention as a promising middle ground between applied mathematics and implicit ML  \cite{JAVADI2025116075}. Unlike purely black-box models, SR searches a space of mathematical expressions to discover parsimonious closed-form relationships that balance interpretability, generalization, and predictive accuracy \cite{udrescu2020ai,cranmer2020discovering}. This makes SR particularly appealing for power system optimization, where tractability and transparency are as critical as predictive performance. By integrating SR into the surrogate modeling for JCC-OPF, one can derive constraints that are not only computationally efficient but also physically meaningful, bridging the gap between theoretical rigor and practical applicability.

\definecolor{masterBlue}{RGB}{27,103,170}
\definecolor{leftGreen}{RGB}{44,160,44}
\definecolor{rightPurple}{RGB}{148,103,189}

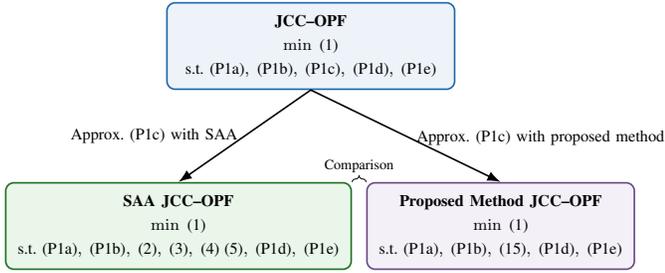
\begin{figure}[!t]
\centering
\resizebox{1\columnwidth}{!}{%
\begin{tikzpicture}[
  >=Latex,
  every node/.style={font=\footnotesize},
  box1/.style={draw=masterBlue, fill=masterBlue!8, rounded corners, thick, align=center, inner xsep=6pt, inner ysep=6pt},
  box2/.style={draw=leftGreen!70!black, fill=leftGreen!10, rounded corners, thick, align=center, inner xsep=6pt, inner ysep=6pt},
  box3/.style={draw=rightPurple!70!black, fill=rightPurple!10, rounded corners, thick, align=center, inner xsep=6pt, inner ysep=6pt},
  arr/.style={-Latex, thick} 
]

\node[box1, minimum width=5.0cm] (master) {
  \textbf{JCC--OPF}\\[1mm]
  $\min\; \eqref{eq:obj}$\\[1mm]
  s.t.\ $\eqref{eq:power_balance},\;\eqref{eq:gen_limits},\;
  \eqref{eq:chance_constraint},\;\eqref{eq:agc_weights},\;\eqref{eq:theta_ref}$
};

\coordinate (splitL) at ($(master.south)+(0,0)$);
\coordinate (splitR) at ($(master.south)+(0,0)$);

\node[box2, minimum width=2cm, below=1.6cm of master, xshift=-2.3cm] (leftbox) {
  \textbf{SAA JCC--OPF}\\[1mm]
  $\min\; \eqref{eq:obj}$\\[1mm]
  s.t.\ $\eqref{eq:power_balance},\;\eqref{eq:gen_limits},\;
  \eqref{eq:active_power_transmission},\;\eqref{eq:big_m_flow},\;
  \eqref{eq:violation_budget}\;\eqref{eq:binary_var},\;
  \eqref{eq:agc_weights},\;\eqref{eq:theta_ref}$
};

\node[anchor=east] at ($(splitL)!0.5!(leftbox.north)$)
  {Approx. \eqref{eq:chance_constraint} with SAA};
\draw[arr] (splitL) -- (leftbox.north);

\node[box3, minimum width=3.8cm, below=1.6cm of master, xshift=3.3cm] (rightbox) {
  \textbf{Proposed Method JCC--OPF}\\[1mm]
  $\min \; \eqref{eq:obj}$ \\[1mm]
  s.t.\ $\eqref{eq:power_balance},\;\eqref{eq:gen_limits},\;
  \eqref{eq:svm},\;\eqref{eq:agc_weights},\;\eqref{eq:theta_ref}$
};

\node[anchor=west] at ($(splitR)!0.52!(rightbox.north)$)
  {Approx. \eqref{eq:chance_constraint} with proposed method};
\draw[arr] (splitR) -- (rightbox.north);

\draw[decorate, decoration={brace, amplitude=3pt}]
  (leftbox.north east) -- (rightbox.north west)
  node[midway, yshift=8pt] {\scriptsize Comparison};

\end{tikzpicture}%
}
\caption{Conceptual comparison between the conventional SAA and the proposed surrogate-based formulation for solving the JCC–OPF problem under renewable uncertainty. The proposed method replaces the scenario-based approximation with a data-driven linear surrogate to enhance tractability and interpretability.
}
\label{fig:jcc-saa-split}
\end{figure}

This paper presents a surrogate approach for JCC-OPF under renewable uncertainty by directly integrating ensemble SVM classifier hyperplane expressions into the optimization problem. In this approach, SVMs yield polyhedral approximations of the feasible region among system variables in interpretable algebraic forms. The resulting surrogate constraints are embedded in the OPF model using Big-M formulations. The proposed methodology is tested on the IEEE 118-bus system under wind and load uncertainty scenarios, with comparative analyses against the conventional scenario approximation approach. The results demonstrate that the surrogate framework achieves a favorable balance between cost optimality, interpretability, and reliability. It represents a promising direction for integrating learning techniques with optimization in renewable-dominated power systems.

\section{Mathematical Framework for JCC-OPF}
This section presents the mathematical formulation of the JCC-OPF problem with affine automatic generation control (AGC) policies. We first introduce the complete optimization model and then provide a detailed exposition of its components, including uncertainty characterization, operational constraints, and the sample average approximation reformulation using a Big-M methodology.

Fig.~\ref{fig:jcc-saa-split} summarizes the reformulation of the JCC–OPF problem under renewable uncertainty. The original JCC–OPF is decomposed into two approaches: the conventional SAA, which replaces probabilistic constraints with scenario-based deterministic ones, and the proposed surrogate-based formulation, which substitutes these constraints with a data-driven linear approximation. This comparison highlights the transition from sampling-based to learning-augmented representations of uncertainty in power system optimization.

\subsection{Problem Formulation}
Consider a power network represented by a directed graph $G(\mathcal{N}, \mathcal{L})$, where $\mathcal{N}$ denotes the set of buses and $\mathcal{L}$ denotes the set of transmission lines. Let $\mathcal{G} \subseteq \mathcal{N}$ represent the subset of buses hosting generation units. The JCC-OPF problem is formulated as \cite{porras2023tight}:

\begin{align}
    \min_{\{p_g, \beta_g, z_s\}} \quad & \sum_{g \in \mathcal{G}} \left( C_{2,g} p_g^2 + C_{1,g} p_g + C_{0,g} + \mathrm{Var}(\Omega) C_{2,g} \beta_g^2 \right) \label{eq:obj}\\
    \text{s.t.} \quad
    & \sum_{g \in \mathcal{G}} p_g + \sum_{n \in \mathcal{N}} \bar{w}_n = \sum_{n \in \mathcal{N}} \bar{d}_n, \tag{P1a} \label{eq:power_balance} \\
    & \underline{p}_g \leq p_g + \beta_g \Omega \leq \overline{p}_g, \quad \forall g \in \mathcal{G}, \tag{P1b} \label{eq:gen_limits} \\
    & \mathbb{P}\Big( |f_l(\omega)| \leq \overline{f}_l, \ \forall l \in \mathcal{L} \Big) \geq 1 - \alpha, \tag{P1c} \label{eq:chance_constraint} \\
    & \sum_{g \in \mathcal{G}} \beta_g = 1, \quad 0 \leq \beta_g \leq 1, \quad \forall g \in \mathcal{G}, \tag{P1d} \label{eq:agc_weights}\\
    & \theta_{ref} = 0, \tag{P1e} \label{eq:theta_ref}
\end{align}
where \eqref{eq:obj} defines the expected total system cost, where $C_{0,g}$, $C_{1,g}$, and $C_{2,g}$ are the quadratic cost coefficients. \eqref{eq:power_balance} enforces power balance. The demand at bus $n \in \mathcal{N}$ is characterized by $d_n = \bar{d}_n + \omega_n$,
where $\bar{d}_n$ represents the nominal demand forecast and $\omega_n$ denotes a zero-mean stochastic forecast error with variance $\sigma_n^2$. The aggregate system-level power imbalance is defined as $\Omega = \sum_{n \in \mathcal{N}} \omega_n$. \eqref{eq:gen_limits} ensures that the limits of the generator capacity are met in all executions. To maintain system balance under uncertainty, each generator uses an affine control policy given by $\tilde{p}_g = p_g + \beta_g \Omega$, where $p_g$ denotes the scheduled base-point generation and $\beta_g \in [0,1]$ represents the AGC participation factor of generator $g$. This linear feedback policy ensures real-time compensation of system imbalances. \eqref{eq:chance_constraint} represents the joint probabilistic transmission constraint, where $\alpha \in (0,1)$ specifies the acceptable risk level and $\overline{f}_l$ denotes the thermal capacity of line $l$. \eqref{eq:agc_weights} enforces AGC participation factor requirements, while \eqref{eq:theta_ref} fixes the reference angle to zero.

\subsection{Sample Average Approximation (SAA) Reformulation}

To address the constraint in \eqref{eq:chance_constraint}, the underlying probabilistic condition is discretized into a finite set of realizations representing sampled uncertainty scenarios. The sample average approximation (SAA) method is then used to approximate the original stochastic constraint by its empirical counterpart, enabling a deterministic and computationally tractable reformulation that can be efficiently solved using standard optimization techniques~\cite{shapiro2021lectures}. 

Using the DC power flow approximation in \eqref{eq:chance_constraint}, the active power flow on transmission line $l \in \mathcal{L}$ under uncertainty scenario $s \in N$ is expressed as 
\begin{align}
f_l(\omega^s) = \sum_{n \in \mathcal{N}} B_{ln} \left( \sum_{g \in \mathcal{G}_n} (p_g + \beta_g \Omega^s) + w_n^s - d_n^s \right),
\label{eq:active_power_transmission}
\end{align}
where $B_{ln}$ denotes the power transfer distribution factor coefficient relating net injection at bus $n$ to flow on line $l$, and $\mathcal{G}$ represents the set of generators located at bus $n$. To render this constraint computationally tractable, we use the SAA method with a Big-M reformulation. Given $N$ independent samples $\{\omega^1, \omega^2, \ldots, \omega^N\}$ drawn from the uncertainty distribution, the JCC is approximated by:

\begin{align}
    -\overline{f}_l - M z_s &\leq f_l(\omega^s) \leq \overline{f}_l + M z_s \quad \forall l \in \mathcal{L}, \ \forall s \in \{1,\dots,N\}, \label{eq:big_m_flow}\\
    \sum_{s=1}^N z_s &\leq \lfloor \alpha N \rfloor, \label{eq:violation_budget}\\
    z_s &\in \{0,1\} \quad \forall s \in \{1,\dots,N\}, \label{eq:binary_var}
\end{align}
where the binary variable $z_s$ indicates whether scenario $s$ is permitted to violate the line flow constraints ($z_s = 1$) or must satisfy them ($z_s = 0$), $M$ is a sufficiently large constant, and $\lfloor \alpha N \rfloor$ represents the maximum number of allowable constraint violations across all samples, thereby approximating the confidence level $1-\alpha$.

This formulation is consistent with state-of-the-art approaches to JCC-OPF in the literature. In particular, it extends the canonical frameworks in \cite{bienstock2014chance,roald2017chance} by incorporating affine AGC policies and using the SAA with a Big-M reformulation for tractability. The resulting model balances fidelity and tractability, enabling systematic enforcement of probabilistic transmission constraints under renewable energy uncertainty.

\section{Explicit Ensemble Support Vector Classification Framework}

This section establishes the theoretical foundation for the ensemble-based classification methodologies used in this work. We begin with the fundamental principles of support vector machines, extend the discussion to ensemble learning paradigms—specifically the bootstrap aggregating (bagging) technique for ensemble construction—and conclude with the formulation of the optimization problem with the ensemble-learned surrogate model.

\subsection{Support Vector Classification for Feasibility Surrogates}

Support vector machines (SVMs) represent a robust class of supervised learning algorithms grounded in statistical learning theory and structural risk minimization \cite{vapnik1995support}. In this work, linear SVMs are used to classify feasible and infeasible operating points of the power system under uncertainty. The training dataset, denoted by $\mathcal{D} = \{(\mathbf{x}_i, y_i)\}_{i=1}^n$, is constructed from samples generated through repeated JCC-OPF simulations, where $\mathbf{x}_i \in \mathbb{R}^d$ represents the system feature vector corresponding to the optimal generator dispatch values ($P_g$), and $y_i \in \{-1, +1\}$ indicates whether the corresponding operating point satisfies all line flow constraints.

The objective of a linear SVM is to find an optimal separating hyperplane that maximizes the margin between feasible ($y_i = +1$) and infeasible ($y_i = -1$) operating regions. The classifier is defined by
\begin{equation}
f(\mathbf{x}) = \mathbf{w}^\top \mathbf{x} + b,
\label{eq:svm_hyperplane_power}
\end{equation}
where $\mathbf{w} \in \mathbb{R}^d$ is the weight vector to the hyperplane and $b \in \mathbb{R}$ is the bias term. The decision function $\text{sgn}(f(\mathbf{x}))$ partitions the input space into feasible and infeasible regions, effectively providing a polyhedral approximation of the system’s true security boundary.

For linearly separable data, the optimal hyperplane can be obtained by maximizing the geometric margin, which is equivalent to solving
\begin{align}
    \min_{\mathbf{w}, b} \quad & \frac{1}{2} \|\mathbf{w}\|^2 \label{eq:svm_obj}\\
    \text{s.t.} \quad & y_i(\mathbf{w}^\top \mathbf{x}_i + b) \geq 1, \quad \forall i \in \{1,\dots,n\}.
    \label{eq:svm_constraint}
\end{align}

However, since the feasible and infeasible regions in power systems may overlap due to uncertainty and nonlinear interactions, a soft-margin SVM is adopted by introducing slack variables $\xi_i \geq 0$ and a regularization parameter $C > 0$:
\begin{align}
    \min_{\mathbf{w}, b, \boldsymbol{\xi}} \quad & \frac{1}{2} \|\mathbf{w}\|^2 + C \sum_{i=1}^n \xi_i \label{eq:soft_svm_obj}\\
    \text{s.t.} \quad & y_i(\mathbf{w}^\top \mathbf{x}_i + b) \geq 1 - \xi_i, \quad \forall i \in \{1,\dots,n\}, \label{eq:soft_svm_constraint1}\\
    & \xi_i \geq 0, \quad \forall i \in \{1,\dots,n\},
    \label{eq:soft_svm_constraint2}
\end{align}
where the parameter $C$ controls the tradeoff between maximizing the separation margin and minimizing misclassifications. In the context of feasibility classification, support vectors correspond to operating points lying near the security boundary—those that most influence the shape of the surrogate feasible region. These hyperplanes collectively serve as linear surrogates of the nonlinear JCCs governing power flow limits.

\subsection{Ensemble Learning for Surrogate Model Construction}

Ensemble learning provides a meta-framework that integrates multiple base classifiers to achieve improved predictive accuracy, robustness, and generalization compared to individual learners \cite{dietterich2000}. In the context of power system optimization, ensemble techniques are particularly valuable for approximating complex nonlinear feasibility boundaries that arise from the interaction between uncertain renewable generation, network topology, and operational constraints. By aggregating multiple weak classifiers, ensemble learning mitigates the effect of individual model bias and variance, resulting in a more reliable surrogate representation of the feasible operating region.

Formally, let $\mathcal{E} = \{h_1, h_2, \ldots, h_M\}$ denote an ensemble comprising $M$ base classifiers, where each $h_m: \mathbb{R}^d \rightarrow \{-1, +1\}$ represents a linear SVM trained on a bootstrap subset of the labeled dataset. The ensemble prediction for a given operating condition $\mathbf{x}$ is obtained via majority voting:
\begin{equation}
    H(\mathbf{x}) = \text{sgn}\left(\sum_{m=1}^M w_m h_m(\mathbf{x})\right),
    \label{eq:ensemble_voting}
\end{equation}
where $w_m \ge 0$ denotes the weight of the $m$-th classifier, typically satisfying $\sum_{m=1}^{M} w_m = 1$. In this study, uniform voting ($w_m = 1/M$) is adopted to ensure equal contribution from all base learners.

The performance of an ensemble surrogate depends on two essential properties: (i) the \textit{accuracy} of each base classifier and (ii) the \textit{diversity} among classifiers. Accuracy ensures that individual SVMs correctly separate feasible and infeasible operating points within their respective bootstrap samples, while diversity guarantees that the classifiers capture different regions of the nonlinear feasibility boundary. This diversity-driven aggregation results in a polyhedral surrogate that provides a more comprehensive and less biased approximation of the joint chance-constrained feasible set \cite{kuncheva2003}.

\subsection{Bootstrap Aggregating for Ensemble SVM Surrogates}

Bootstrap aggregating, or bagging \cite{breiman1996}, is an ensemble learning technique that mitigates variance and enhances model stability by training multiple base classifiers on randomly resampled subsets of the training data. In the context of power system optimization, bagging enables the construction of a robust surrogate that captures the complex and nonlinear feasible boundary of the JCC-OPF problem under renewable uncertainty. By generating diverse SVM classifiers trained on different bootstrap samples, the ensemble achieves a more generalized and less overfitted approximation of the underlying feasible set.

\subsubsection{Bootstrap Sampling}

Given a labeled dataset $\mathcal{D} = \{(\mathbf{x}_i, y_i)\}_{i=1}^n$, derived from JCC-OPF simulations under varying wind and load scenarios, a bootstrap sample $\mathcal{D}_m^*$ is constructed by sampling $n$ instances from $\mathcal{D}$ with replacement. Each bootstrap dataset typically contains about $63.2\%$ unique samples from the original set, while the remaining $36.8\%$—known as out-of-bag samples—can be used for validation and unbiased performance estimation. This stochastic resampling process ensures that each base classifier in the ensemble learns slightly different decision boundaries, promoting diversity among the weak SVM models and reducing correlation in their predictions.

\subsubsection{Bagging SVM Algorithm}

The bagging-based ensemble SVM surrogate used in this study follows the steps below:

\begin{enumerate}
\item For each bootstrap index $m = 1, \ldots, M$:
\begin{enumerate}
\item Generate bootstrap sample $\mathcal{D}_m^* = \{(\mathbf{x}_i^*, y_i^*)\}_{i=1}^n$ by sampling with replacement from $\mathcal{D}$.
\item Train a linear SVM classifier $h_m$ on $\mathcal{D}_m^*$ using the soft-margin formulation (\eqref{eq:soft_svm_obj}--\eqref{eq:soft_svm_constraint2}) to obtain its separating hyperplane.
\end{enumerate}
\item Aggregate the predictions of the $M$ trained SVMs using majority voting in \eqref{eq:ensemble_voting}.
\end{enumerate}

The resulting ensemble classifier is expressed as:
\begin{equation}
    H_{\text{Bagging-SVM}}(\mathbf{x}) = \text{sgn}\left(\frac{1}{M}\sum_{m=1}^M h_m(\mathbf{x})\right),
    \label{eq:bagged_svm}
\end{equation}
where each $h_m(\mathbf{x}) = \mathbf{w}_m^\top \mathbf{x} + b_m$ denotes the decision of the $m$-th SVM trained on a bootstrap dataset.

The theoretical advantage of bagging lies in its ability to reduce variance without substantially increasing bias. Assuming $M$ base classifiers with individual error rate $\epsilon$ and pairwise correlation $\rho$, the expected error rate of the ensemble can be approximated as:
\begin{equation}
\epsilon_{\text{bagging}} \approx \rho \epsilon + \frac{1-\rho}{M}\epsilon,
\label{eq:bagging_error_power}
\end{equation}
indicating that the ensemble error decreases as $M$ grows, converging toward $\rho \epsilon$ as $M \to \infty$. In power system applications, this property is particularly beneficial for maintaining reliability in surrogate representations of probabilistic constraints, where uncertainty-induced nonlinearities can otherwise lead to inconsistent classification of feasible and infeasible operating points.

\subsection{Integration of the Surrogate Model}

Consider the JCC-OPF problem introduced in Section~II. In this work, a trained surrogate model based on ensemble linear SVM classifier is used to approximate the feasibility boundary of the original design space. The surrogate replaces the probabilistic constraint in \eqref{eq:chance_constraint}. To account for a predefined risk level, a Big--M reformulation is applied, which permits controlled violations of the surrogate feasibility region. Since the original JCC is approximated using $N$ sampled scenarios, the surrogate formulation adopts the same scenario-based structure for consistency. The resulting SVM-based surrogate JCC-OPF with Big--M integration is formulated as follows:

Let $\{(w_m,b_m)\}_{m=1}^M$ be the ensemble SVM hyperplanes (trained offline) acting on the redispatched outputs ($P_{g}$). For each scenario,
\begin{equation}
\frac{1}{M}\sum_{m=1}^M(w_m^\top P_{g,s} + b_m) + M_{\mathrm{svm}} z_s \;\ge\; 0,
\qquad \forall s \in \{1,\dots,N\}, \label{eq:svm}
\end{equation}
with a sufficiently large scalar $M_{\mathrm{svm}}$. If $z_s=0$, all hyperplanes are enforced (classified safe); if $z_s=1$, the surrogate is relaxed for scenario $s$. The optimization problem incorporating the ensemble linear SVM surrogate can be expressed in the following compact form:  
\begin{align}
    \min \; & \eqref{eq:obj}  \nonumber \\
    \text{s.t.} \; & \eqref{eq:power_balance}, \; \eqref{eq:gen_limits}, \; \eqref{eq:svm}, \; \eqref{eq:agc_weights}, \; \eqref{eq:theta_ref}.
\end{align}

\definecolor{cA}{RGB}{222,235,247}
\definecolor{cB}{RGB}{254,230,206} 
\definecolor{cC}{RGB}{204,235,197} 
\definecolor{cD}{RGB}{233,222,238}

\begin{figure}[t]
\centering
\begin{tikzpicture}[
  >=Latex,
  node distance=6mm and 6mm,
  every node/.style={font=\footnotesize},
  box/.style={
    draw=black!55, rounded corners=2pt, align=center, line width=0.3pt,
    minimum height=10mm, inner sep=6pt
  },
  arrow/.style={-Latex, line width=0.7pt, draw=black!60}
]

\node[box, fill=cA, minimum width=0.60\columnwidth] (A) {%
\textbf{Generate and label data}\\
Generate data from SAA JCC-OPF \\ under wind/load uncertainty and label samples (\(P_g, y)\)
};

\node[box, fill=cB, minimum width=0.60\columnwidth, below=4mm of A] (B) {%
\textbf{Train ensemble SVM surrogate}\\
Train $M$ linear SVMs using bagging \\ Aggregate half–spaces into a polyhedral surrogate
};

\node[box, fill=cD, minimum width=0.60\columnwidth, below=4mm of B] (C) {%
\textbf{Embed surrogate model in JCC-OPF}\\
Replace line flow limits in JCC-OPF with ensemble SVM surrogate \\ Big-M enables risk–budgeted relaxation (same $\alpha$ budget)
};

\node[box, fill=cC, minimum width=0.60\columnwidth, below=4mm of C] (D) {%
\textbf{Solve surrogate JCC-OPF model}\\
Solve surrogate JCC-OPF model and compare with SAA JCC-OPF
};

\draw[arrow] (A) -- (B);
\draw[arrow] (B) -- (C);
\draw[arrow] (C) -- (D);

\end{tikzpicture}
\caption{Overall workflow of the proposed ensemble SVM surrogate approach. Data is generated and labeled under wind/load uncertainty using SAA JCC-OPF. An ensemble of linear SVMs is trained via bagging to construct a polyhedral surrogate, which is then embedded into the JCC-OPF formulation through Big–M relaxation. The surrogate JCC-OPF is subsequently solved and benchmarked against the baseline formulation.}
\label{fig:bigpicture_pipeline}
\end{figure}
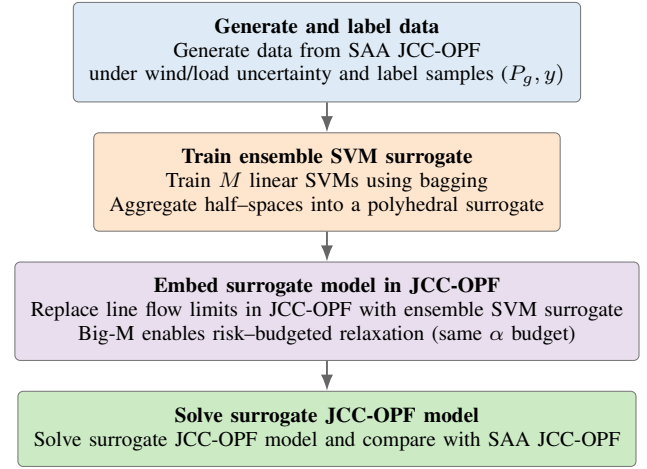

Fig.~\ref{fig:bigpicture_pipeline} illustrates the overall workflow of the proposed surrogate-based framework for solving the JCC-OPF problem under renewable uncertainty. In the first stage, data is generated by repeatedly solving the scenario-based approximation of JCC-OPF with load and wind perturbations, and each solution is labeled according to its feasibility. In the second stage, an ensemble of linear SVMs is trained using bagging, where multiple bootstrap samples are used to construct several classifiers that are aggregated into a polyhedral surrogate of the feasibility boundary. In the third stage, the surrogate is embedded into the JCC-OPF formulation by replacing the original line-flow chance constraints with SVM-derived half-spaces, while using a Big–M relaxation to allow controlled risk according to the specified violation budget. The surrogate JCC-OPF is solved, providing an approximation of the baseline problem that preserves key reliability and tractability properties.

\section{Numerical results \& discussion}
The simulation results are generated using the IEEE 118-bus test system. The JCC-OPF problem is implemented in \textsc{Matlab} using \textsc{CPLEX} with the \textsc{MATPOWER} and \textsc{YALMIP} toolboxes. The ML tasks are performed in \textsc{Python}. All computations are performed on a personal computer equipped with an Intel(R) Core(TM) Ultra 9 285K CPU @ 3.70\,GHz, 32\,GB of RAM, and a 64-bit operating system. The results obtained from the proposed method are compared with those derived from SAA to provide deeper insights into performance. 

\subsection{Dataset Generation for Ensemble SVM Training}

To train the surrogate classifier, a labeled dataset was generated using a scenario-based approximation of the JCC-OPF with AGC participation. Each data point corresponds to a sampled wind profile, an associated OPF solution, and an ex-post feasibility label.

To model load uncertainty, each bus demand was perturbed according to
\begin{equation}
    P_{d,s} = P_d^{\mathrm{mean}} + \sigma_d \odot \xi_d, \qquad 
    \xi_d \sim \mathcal{N}(0, I),
\end{equation}
where $P_d^{\mathrm{mean}}$ is mean active power demand (load), $\sigma_d = 0.03\,P_d^{\mathrm{mean}}$ denotes the standard deviation of load fluctuations, $\odot$ is the element-wise (Hadamard) multiplication, and $\xi_d$ is a random perturbation vector. 

Wind uncertainty was introduced at selected generator buses $\mathcal{N}_w = \{9,31,54,90,100\}$ with nominal mean $\mu_w$ and standard deviation $\sigma_w$. For each run, the nominal wind statistics were perturbed to create diverse operating conditions. Specifically, the mean and standard deviation of wind generation at the selected buses were randomly varied within realistic bounds, ensuring that the resulting values remained physically meaningful and representative of system uncertainty. These perturbed parameters were then used to define the probability distributions from which wind injection scenarios were sampled.

For each run, $r=1,\dots,n_{\text{runs}}$, $N=100$ scenarios are sampled:
\begin{equation}
    P_{w,s} = \mu_{w}^{(r)} + \sigma_{w}^{(r)} \odot \xi_{w,s}, \quad
    P_{d,s} = P_{d}^{\mathrm{mean}} + \sigma_d \odot \xi_{d,s},
\end{equation}
with $\xi_{w,s}, \xi_{d,s}\sim \mathcal{N}(0,I)$. Here, $P_{w,s}$ denotes the wind power injection vector, $\mu_{w}^{(r)}$ is the mean wind generation vector, $\sigma_{w}^{(r)}$ is the standard deviation vector, and $\xi_{w,s}$ are independent Gaussian random variables. For each sampled wind profile, two JCC-OPF problems were solved with violation probability levels $\alpha \in \{0,0.05\}$, yielding the optimized generator dispatch and AGC participation factors.

After generating datasets corresponding to each violation probability level, the line flow constraints were evaluated against their respective limits. Ex-post feasibility was further validated using an independent Monte Carlo test set to ensure accurate labeling. Each scenario was labeled as \textit{feasible ($y=+1$)} if no violations occurred, and \textit{infeasible ($y=-1$)} if at least one line constraint was breached. To mitigate class imbalance, feasible and infeasible samples were selectively resampled until both classes were represented in approximately equal proportion.

Each feature vector included only the optimized generator dispatch levels ($P_g$), which can capture both the uncertainty inputs and the system’s operational response. The dataset was then partitioned into 75\% for training and 25\% for testing. The training portion was used to construct an ensemble surrogate model by applying the bagging technique. Specifically, eight weak linear SVM classifiers were trained, and their outputs were aggregated to form an ensemble SVM. This ensemble size was found to offer a favorable balance between variance reduction and computational cost. The resulting hyperplanes were subsequently used as surrogate representations of the line flow limits under JCCs. All random sampling procedures were performed with fixed seeds to ensure reproducibility of the dataset and training process.
 
\subsection{Results for IEEE 118-bus system}
By replacing the surrogate model with the line flow limits in the JCC-OPF formulation and evaluating it on 15 test samples, the resulting operating costs were compared with those obtained from the SAA of JCC-OPF using the Big--M method. The comparative results are presented in Table~\ref{comp}. Across all 15 test samples, the surrogate consistently produces solutions whose costs closely track the original JCC-OPF benchmark. Although the surrogate method incurs a slight cost penalty relative to the exact formulation, the gap remains small, with relative deviations of less than $0.05\%$ on average.

\begin{table}[t]
\centering
\caption{Cost comparison between JCC-OPF and Ensemble SVM across 15 test samples.}
\label{tab:cost_comparison}
\begin{tabular}{ccccc}
\toprule
\textbf{Sample} & \textbf{JCC-OPF(\$)} & \textbf{Ensemble SVM(\$)} & \textbf{$\Delta$Cost(\$)} & \textbf{$\Delta$Cost(\%)} \\
\midrule
1 & 90772 & 90814 & 42 & 0.046\\
2 & 90781 & 90821 & 40 & 0.044\\
3 & 90795 & 90834 & 39 & 0.044\\
4 & 90836 & 90873 & 37 & 0.040\\
5 & 90862 & 90898 & 36 & 0.040\\
6 & 90884 & 90921 & 37 & 0.041\\
7 & 90946 & 90979 & 33 & 0.036\\
8 & 90992 & 91022 & 30 & 0.033\\
9 & 91046 & 91073 & 27 & 0.030\\
10 & 91055 & 91082 & 27 & 0.029\\
11 & 91062 & 91087 & 25 & 0.027\\
12 & 91076 & 91096 & 20 & 0.023\\
13 & 91084 & 91109 & 25 & 0.027\\
14 & 91106 & 91130 & 24 & 0.026\\
15 & 91301 & 91315 & 14 & 0.015\\
\bottomrule
\end{tabular}
\label{comp}
\end{table}

The results demonstrate that the ensemble SVM can reliably approximate the feasibility boundary imposed by JCCs while preserving near-optimal cost efficiency. This tradeoff highlights the strength of the surrogate approach: by replacing computationally expensive probabilistic line flow constraints with a polyhedral approximation learned from data, the method achieves scalability and tractability at the expense of only marginal increases in operating cost.

\begin{figure}[t]
    \centering
    \includegraphics[width=0.7\linewidth]{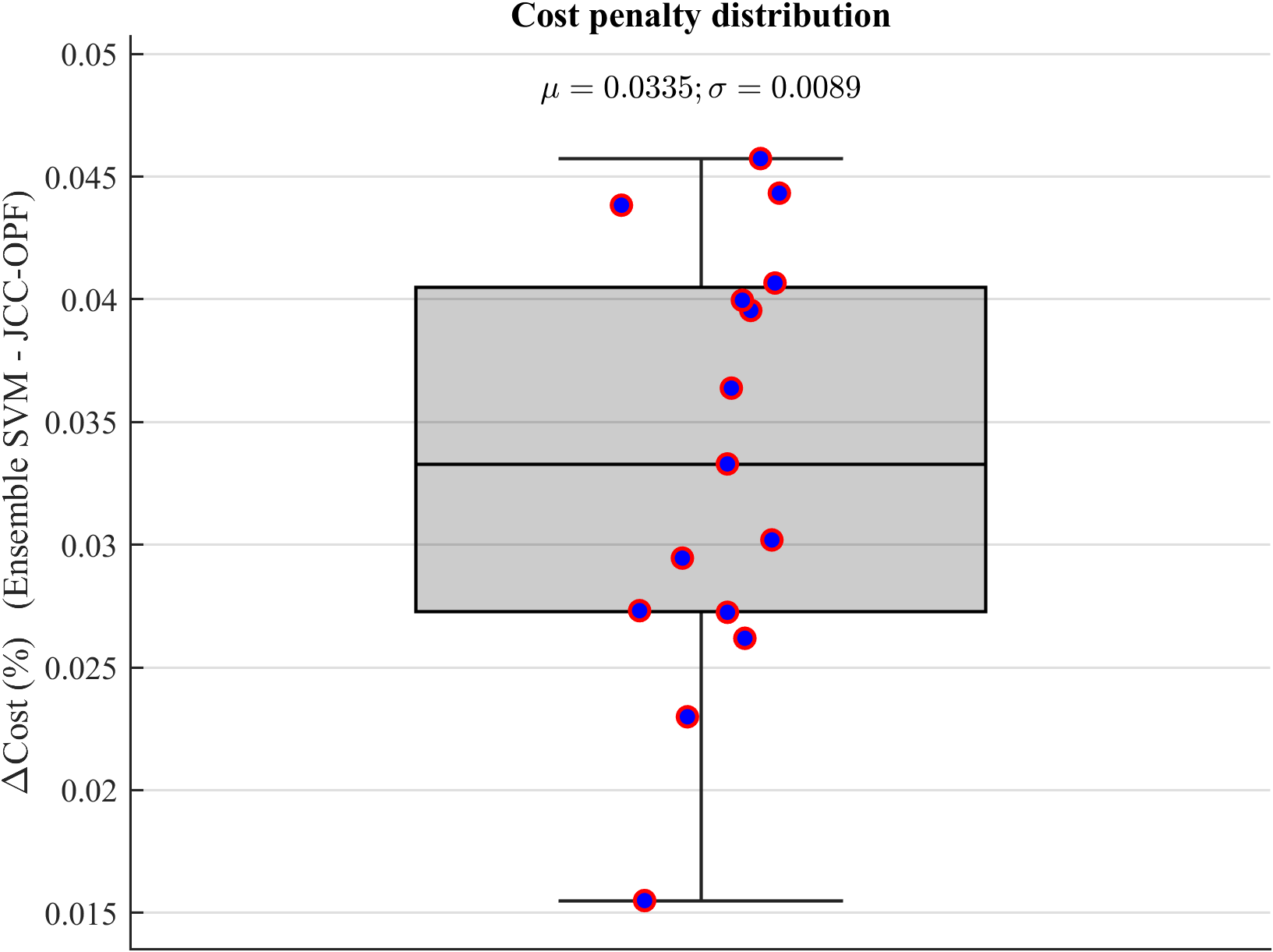}
    \caption{Distribution of the relative cost penalty of the ensemble SVM surrogate with respect to the baseline JCC-OPF. The boxplot summarizes results across 15 test samples, with individual sample points shown in red. The mean and standard deviation are $\mu=0.0335\%$ and $\sigma=0.0089\%$, respectively.}
    \label{cost_delta}
\end{figure}

Fig.~\ref{cost_delta} illustrates the distribution of the relative cost penalty incurred when using the ensemble SVM surrogate in place of the exact JCC-OPF formulation. Across the 15 test samples, the surrogate consistently yielded slightly higher costs compared to the baseline, with a mean penalty of $0.0335\%$ and a standard deviation of $0.0089\%$. In absolute terms, this corresponds to an average increase of approximately \$30, with a variability of \$8, ranging from a minimum of \$14 to a maximum of \$42.

The interquartile range remains narrow, indicating that the additional cost burden is stable across different operating conditions. The results confirm that the surrogate model achieves a high degree of fidelity in approximating the feasibility region of the JCC-OPF. Although minor cost increases are observed, these are negligible in practical terms, especially given the significant computational benefits of replacing probabilistic line flow constraints with data-driven surrogate hyperplanes. Thus, the ensemble SVM provides an effective tradeoff between efficiency and operational optimality.

\begin{figure}[t]
    \centering
    \includegraphics[width=0.7\linewidth]{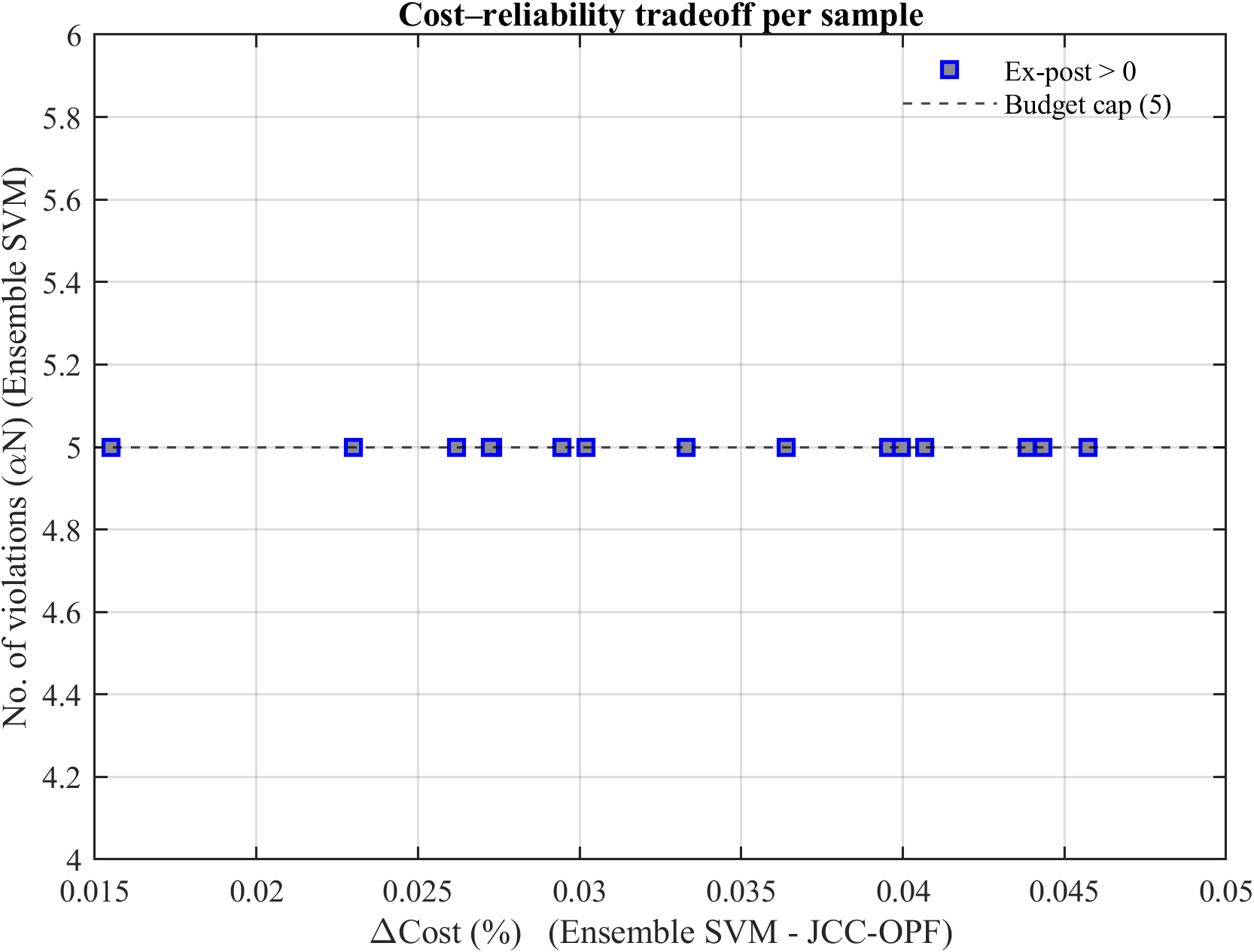}
    \caption{Cost–reliability tradeoff for the ensemble SVM surrogate across 15 test samples. The $x$-axis shows the relative cost penalty with respect to JCC-OPF, while the $y$-axis indicates the number of ex-post violations. All points lie on the violation budget cap ($\alpha N = 5$), demonstrating that the surrogate enforces the prescribed reliability level while incurring only marginal cost increases.}
    \label{tradeoff}
\end{figure}

Fig.~\ref{tradeoff} presents the tradeoff between cost and reliability when applying the ensemble SVM surrogate. The $x$-axis reports the relative cost penalty with respect to the exact JCC-OPF, while the $y$-axis captures the number of ex-post violations observed across 100 scenarios per sample.

The results indicate that for all 15 test cases, the surrogate consistently operated at the violation budget cap of $\alpha N = 5$, matching the prescribed JCC level of $\alpha = 0.05$. This confirms that the surrogate does not compromise reliability, as no test case exceeded the allowable violation threshold. On the other hand, the cost penalty remains small, ranging from approximately $0.015\%$ to $0.045\%$ relative to the baseline, aligning with the distribution reported in Fig.~\ref{cost_delta}.

Overall, these results demonstrate that the ensemble SVM surrogate achieves an effective cost–reliability balance, strictly enforcing the probabilistic security level while introducing only negligible increases in operational cost. This demonstrates the surrogate’s robustness and practicality for large-scale chance-constrained OPF problems.

\begin{figure}[t]
    \centering
    \includegraphics[width=0.7\linewidth]{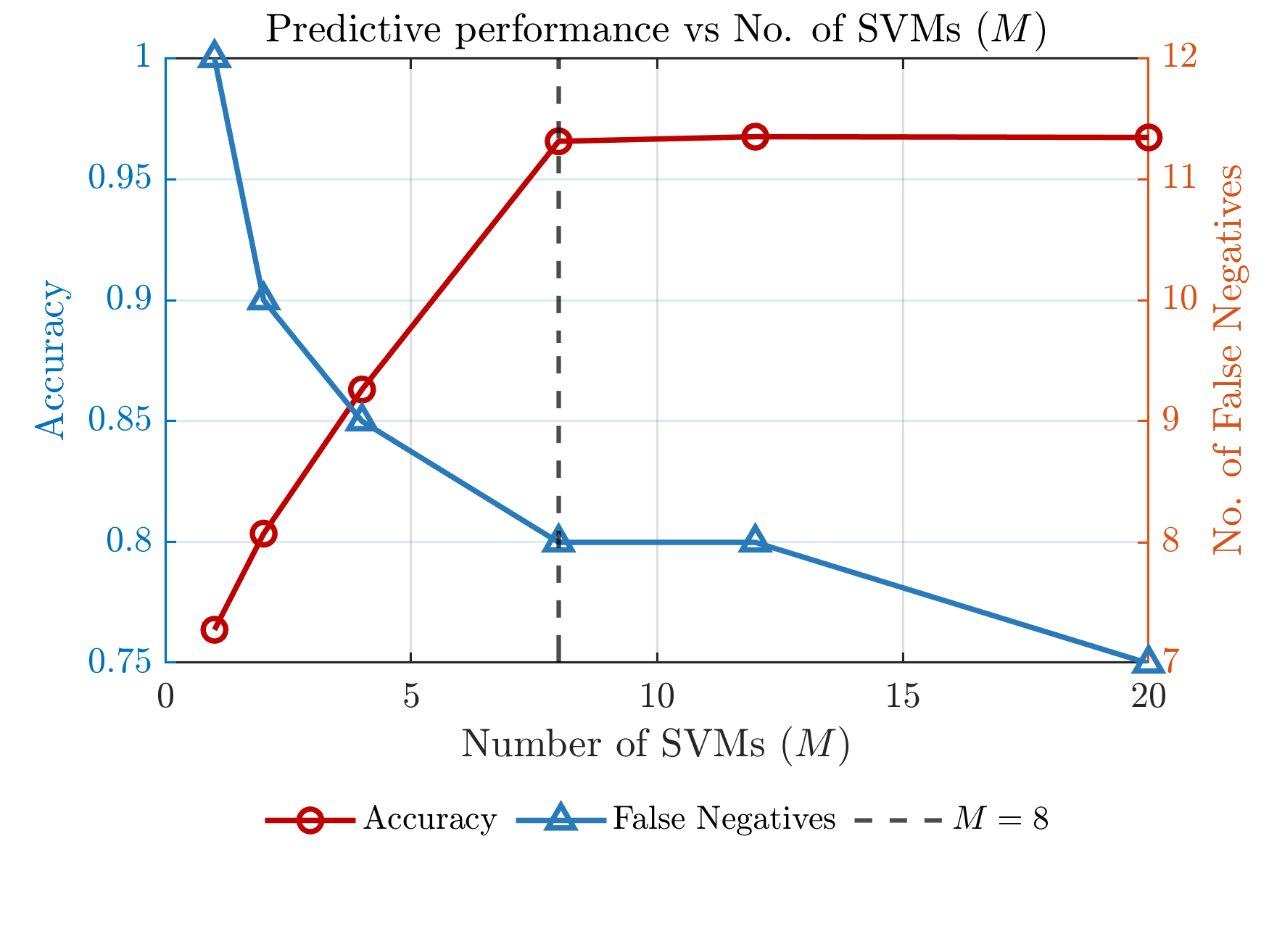}
    \caption{Predictive performance of the ensemble surrogate as a function of the number of weak linear SVMs. Accuracy (red line) improves with ensemble size, while the number of false negatives (blue line) decreases. Beyond $M=8$, gains in predictive performance are marginal, suggesting that eight classifiers offer a favorable balance between accuracy and reliability.}
    \label{performance_B}
\end{figure}

Fig.~\ref{performance_B} illustrates how the predictive accuracy and false negative rate of the ensemble surrogate evolve as the number of weak linear SVMs increases. The results show a clear improvement in both metrics as the ensemble size grows from 1 to 8 classifiers, indicating that bagging effectively reduces variance and improves generalization.

However, beyond $M=8$, the incremental gains in accuracy are negligible, and the false negative rate stabilizes around its minimum value. At the same time, increasing the ensemble size further would proportionally increase the number of surrogate hyperplanes that must be embedded in the OPF model, leading to a higher computational burden without a significant improvement in predictive performance.

Thus, selecting $M=8$ weak linear SVMs provides a strong tradeoff: it achieves near-optimal classification accuracy, reduces false negatives to acceptable levels, and avoids unnecessary complexity in the optimization stage. This balance makes the ensemble with eight learners a practical and scalable surrogate choice for JCC-OPF applications.

\subsection{Discussion}

The results demonstrate that the proposed ensemble SVM surrogate is an effective and scalable approximation of the JCC-OPF problem. By replacing probabilistic line flow constraints with polyhedral surrogates learned from data, the framework achieves a favorable balance between cost optimality and reliability. Across 15 samples from the test set, the surrogate consistently produced near-optimal operating costs, with an average cost penalty of only $0.0335\%$ (approximately \$30) relative to the SAA of JCC-OPF with Big--M reformulation. The distribution of cost penalties further showed a narrow spread, confirming that the surrogate introduces only negligible economic deviations under varying realizations of uncertainty.

Reliability assessment through ex-post validation confirmed that the surrogate respects the prescribed violation budget. The cost--reliability tradeoff analysis indicated that the method achieves robust performance without exceeding the allowable violation rate of $\alpha=0.05$, while incurring only marginal cost increases. Importantly, ensemble learning played a key role in mitigating the risk of misclassification. The tradeoff analysis of ensemble size showed that predictive performance improves with the number of weak SVM classifiers, but stabilizes beyond eight learners. Thus, $M=8$ was identified as a suitable ensemble size, striking a balance between classification accuracy and false negative reduction.

Overall, these findings validate the feasibility of integrating ensemble-based surrogates into large-scale power system optimization under uncertainty. The surrogate approach maintains interpretability, provides transparent hyperplane representations of feasibility boundaries, and scales efficiently compared to conventional scenario approximation methods.

\section{Conclusion}
This paper presented an ensemble SVM surrogate framework for the JCC-OPF problem under renewable uncertainty. By embedding multiple linear SVM hyperplanes into the optimization model, the surrogate approximates the feasible region defined by probabilistic line flow constraints with high fidelity. Comparative simulations on the IEEE 118-bus system revealed that the surrogate consistently produced solutions that closely tracked those of the SAA with Big--M method, incurring only negligible cost increases while satisfying the prescribed violation budget. Ensemble learning was shown to mitigate misclassification risk, with eight weak classifiers providing a suitable balance between predictive performance and computational complexity. Overall, the results highlight the effectiveness of ensemble surrogates in reconciling the trade-off between tractability, interpretability, and reliability in chance-constrained power system optimization. 

\section{Future Work}


Future research can extend this work by incorporating nonlinear kernels into the surrogate model to better capture complex dynamics, applying the approach to AC-OPF and multi-period settings, and integrating advanced uncertainty representations such as distributionally robust or heavy-tailed models. Additionally, embedding the surrogate in real-time decision-making could enhance practical applicability. Overall, the framework provides a scalable and interpretable basis for advancing risk-aware power system operation.

\bibliographystyle{IEEEtran}
\bibliography{refs}

\end{document}